\documentclass[10pt]{cmslatex}

\usepackage{subfig,latexsym,amsmath,amssymb, amsfonts,enumerate,
 epsfig, graphics,times,color}

\newtheorem{remark}[theorem]{Remark}

\numberwithin{figure}{section}

\newcommand{\ccdot}{\;\cdot\;}
\newcommand{\eqdef}{\!\overset{\mbox{\tiny def}}{=}\!}

\newcommand{\N}{\mathbb{N}}
\newcommand{\one}{\mathbf{1}}
\newcommand{\R}{\mathbb{R}}
\newcommand{\C}{C}
\newcommand{\ym}{y^*}
\newcommand{\eqdist}{\!\overset{\mathcal{D}}{=}\!}
\newcommand{\algNam}{Weak Trapezoidal } 
\newcommand{\algNamFull}{Weak $\theta$-Midpoint Trapezoidal } 
\def\E{\mathbb{E}}
\def\PP{\mathbb{P}}

\def\Ft{\mathcal{F}_{t}}
\def\Ftt{\mathcal{F}_{\theta h}}

\title{A weak trapezoidal method for a class of stochastic
  \\differential equations}

\author{David F. Anderson$^{1}$ and Jonathan C. Mattingly$^{2}$}

\begin{document}

\maketitle

\footnotetext[1]{Department of Mathematics, University of
  Wisconsin-Madison, Madison, Wi. 53706, anderson@math.wisc.edu}

\footnotetext[2]{Department of Mathematics, Center of Nonlinear and
  Complex systems, Center for Theoretical and Mathematical Science,
  and Department of Statistical Science, Duke University, Durham,
  N.C. 27708, jonm@math.duke.edu}

\begin{abstract}
  We present a numerical method for the approximation of solutions for
  the class of stochastic differential equations driven by Brownian
  motions which induce stochastic variation in fixed directions.  This
  class of equations arises naturally in the study of population
  processes and chemical reaction kinetics.  We show that the method
  constructs paths that are second order accurate in the weak sense.
  The method is simpler than many second order methods in that it
  neither requires the construction of iterated It\^o integrals nor
  the evaluation of any derivatives. The method consists of two steps.
  In the first an explicit Euler step is used to take a fractional
  step.  The resulting fractional point is then combined with the initial point
  to obtain a higher order, trapezoidal like, approximation.  The
  higher order of accuracy stems from the fact that both the drift and
  the quadratic variation of the underlying SDE are approximated to
  second order.

\end{abstract}

\section{Introduction}
\label{sec:intro}
We consider the problem of constructing accurate approximations on bounded time intervals to
solutions of the following family of stochastic differential equations (SDEs)
\begin{equation}
 \begin{aligned}
 dX(t) &= b(X(t))dt + \sum_{k = 1}^M \sigma_k(X(t)) \;\nu_k \ dW_k(t),\\
 X(0)&= x \in \R^d
 \end{aligned}
  \label{eq:main}
\end{equation}
where $\ b\colon\R^d \to \R^d$,
$\sigma_k\colon \R^d \to \R_{\ge 0}$, $\nu_k \in \R^d$, and $W_k(t)$
are one-dimensional Wiener processes.  Thus, randomness is entering
the system in fixed directions $\nu_k$, but at variable rates
$\sigma_k(X(t))$.  Precise regularity conditions on the coefficients
will be presented with our main results in Section \ref{sec:method}.

The algorithm developed in this paper is a trapezoidal-type method and
consists of two steps; in the first an explicit Euler step is used to
take a fractional step and in the second the resulting fractional point is used
in  combination with the initial point to obtain a higher
order, trapezoidal like, approximation.  We will prove that the method
developed is second order accurate in the weak sense.  Because the
method developed here produces single paths, it is natural to allow
variable step-sizes; this is in contrast to Richardson extrapolation
techniques (\cite{talay1990}). Finally, it is important to note that
while the method presented in this paper is applicable to only a
sub-class of SDEs, that sub-class does include systems whose diffusion
terms do not commute, which is a classical simplifying assumption to
obtain higher order methods (See
\cite{KloedenPlaten92,MilsteinTretyakov05}).

The method we propose is in some sense similar to the classical
predictor-corrector. There have already been a number of such methods
proposed in the stochastic context to produce higher-order methods
(see \cite{PardouxTalay85,Platen95,BurrageTian02}). In a general way,
all of these methods require the simulation of iterated It\^o
integrals and sometimes need derivatives of the diffusion terms. If
one only cares about weak accuracy, it is possible to use random
variables which make these calculations easier and computationally
cheaper. That being said, the complexity and cost of such calculations
is one of the main impediments to their wider use. By assuming a
certain structure for \eqref{eq:main}, we are able to develop a numerical
method which we hope is more easily applied and implemented.

Though a specific structure of \eqref{eq:main} is assumed, it is a
structure which arises naturally in a number of settings.  For
example, our method will be applicable whenever $d = 1$.  Also, we
note that diffusion approximations to continuous time Markov chain
models of population processes, including (bio)chemical processes,
satisfy \eqref{eq:main}.  As stochastic models of biochemical reaction
systems, and, in particular, gene regulatory systems, are becoming
more prevalent in the science literature, developing algorithms that
utilize the specific structure of such models has increased importance
(\cite{Anderson2007b, AndGangKurtz2009}). Furthermore, in
Section~\ref{uniformE}, we quote a result from the literature which
states that any system with uniformly elliptic diffusion can be put in the form of
\eqref{eq:main} without changing its distribution.

The topic of this paper is a method that produces a weak approximation
rather than a strong approximation in that the approximate trajectory
is produced without reference to an underlining Wiener process
trajectory. We see this as an advantage. Except for applications such
as filtering or certain problems of collective motion for stochastic
flows, one is usually simply interested in generating an accurate draw
from the distribution on $C([0,T], \R^d)$ induced by \eqref{eq:main}.
This is different than accurately reproducing the It\^o map $W \mapsto
X(t,W)$ implied by \eqref{eq:main}. The second is referred to as
strong approximation. In our opinion such approximations are usually
unnecessary and lead to a concept of accuracy which is unnecessarily
restrictive. In \cite{GainesLyons97}, it is discussed that without
accurately estimating second order It\^o integrals one cannot produce
a strong method of order greater then 1/2. If the vector fields
commute, then this restriction does not apply and higher order strong
methods are possible. While the term ``strong approximation'' is quite
specific, the term ``weak approximation'' is used for a number of
concepts. Here we mean that the joint distribution of the numerical
method at a fixed number of time points converges to the true marginal
distribution as the numerical grid converges to zero. If this error
goes to zero as the numerical mesh size to the   power $p$ in some
norm on measure then we say the method is of order $p$.  This should
be contrasted with talking about the rate at which a given function of
the path converges.

The outline of the paper is as follows.  In Section~\ref{sec:method}
we present our algorithm together with our main results concerning its
weak error properties.  In Section~\ref{sec:why} we give the intuition
as to why the method should work.  In Section~\ref{sec:proof} we give
the delayed proof of the local error estimates for the method which
were stated in Section~\ref{sec:method}.  In Section
\ref{sec:examples} we provide examples illustrating the performance of
the proposed algorithm.  In Section~\ref{sec:theta} we discuss the
effect of varying the size of the first fractional step of the
algorithm. In Section~\ref{sec:richardson} we compare one step of the
algorithm to one step in a Richardson extrapolation type algorithm. In
Section~\ref{uniformE} we show how, at least theoretically, the method
can be applied to any uniformly elliptic SDE.  Finally, an appendix
contains a tedious calculation needed in Section~\ref{sec:proof}.

\section{The numerical method and main results}
\label{sec:method}
Throughout the paper, we let $X(t)$ denote the solution to
\eqref{eq:main} and $Y_{i}$ denote the computed approximation at the
time $t_i$ for the time discretization $0 = t_0 < t_1< \cdots$.  We
begin both from the same initial condition, namely $X(0) = Y_0 = x_0$.
Let $ \big\{ \eta_{1k}^{(i)}, \eta_{2k}^{(i)}\colon k
\in\{1,\ldots,M\}, i \in \N\big\}$ be a collection of mutually
independent Gaussian random variables with mean zero and variance one.
It is notationally convenient to define $[x]^+= x \vee 0 =
\max\{x,0\}$.

We propose the following algorithm to approximate the solutions of
\eqref{eq:main}.

\vspace{1ex}

\noindent\textsc{Algorithm}. ( \algNamFull)
\textit{Fixing a \, $\theta \in (0,1)$, we define
 \begin{align}
   \alpha_1 \eqdef \frac{1}{2}\frac{1}{\theta(1 -
     \theta)} 
   \qquad\text{and}\quad \alpha_2 \eqdef
   \frac{1}{2}\frac{(1-\theta)^2 + \theta^2}{\theta(1 - \theta)}\,.
 \end{align}
 Next fixing a discretization step $h$, for each $i \in
 \{1,2,3,\dots\}$ we repeat the following steps in which we first
 compute a $\theta$-midpoint $y^*$ and then the new value $Y_i$:
 \begin{enumerate}
 \item[Step 1.] $\displaystyle \ym = Y_{i-1} + b(Y_{i-1}) \theta h +
   \sum_{k=1}^M \sigma_k( Y_{i-1} ) \;\nu_k \; \eta_{1k}^{(i)}
   \sqrt{\theta h}$
 \item [Step 2.] $\displaystyle Y_i = \ym + (\alpha_1 b(\ym) -
   \alpha_2 b(Y_{i-1}))(1-\theta) h + \sum_{k=1}^M
   \sqrt{\big[\alpha_1 \sigma^2_k(\ym) -
     \alpha_2\sigma^2_k(Y_{i-1})\big]^+} \;\nu_k\; \eta_{2k}^{(i)}
   \sqrt{(1-\theta)h}$.
 \end{enumerate}
}
  
\vspace{1ex}
  
\begin{remark}
  Notice that on the $i$th-step $\ym$ is the standard Euler
  approximation to $X(\theta h + (i-1)h)$ starting from $Y_{i-1}$ at
  time $(i-1)h$ \cite{Milstein95}.
\end{remark}
\begin{remark}
  Notice that for all $\theta \in (0,1)$ one has $\alpha_1 > \alpha_2$
  and $\alpha_1-\alpha_2=1$.  It is reasonable to ask which $\theta$
  is best.  Notice that when $\theta = 1/2$ both $\alpha_1$ and
  $\alpha_2$ are minimized with values $\alpha_1 = 2$ and $\alpha_2 =
  1$.  This likely has positive stability implications.  From the
  point of view of accuracy $\theta=1/2$ also seems like a reasonable
  choice as it provides a central point for building a balanced
  trapezoidal approximation, as will be explained in
  Section~\ref{sec:why}.  Further, picking a $\theta$ close to 1 or 0
  increases the likelihood that the term $[\alpha_1 \sigma^2_k(\ym) -
  \alpha_2\sigma^2_k(Y_{i-1})]^+$ will be zero, which will lower the
  accuracy of the method. If instability due to stiffness is a
  concern, one might consider a $\theta$ closer to one as that would
  likely give better stability properties being closer to an implicit
  method. In general, $\theta=1/2$ seems like a reasonable compromise,
  though this question requires further investigation and will be
  briefly revisited in Section \ref{sec:theta}.
  \label{rem:theta}
\end{remark}

For simplicity, we will restrict ourselves to the case when $b$ and
the $\sigma_k$ are in $\C^6(\R^d)$, the space of bounded functions
whose first through sixth derivatives are continuous and bounded. In
general, we will denote by $\C^k(\R^d)$ the space of bounded,
continuous functions whose first $k$ derivatives are bounded and
continuous.  For $f \in \C^k(\R^d)$, we define the standard norm
\begin{align*}
  \|f\|_k = \sup \big\{ |f(x)|, |\partial_\alpha f(x)| : x \in \R^d,
  \alpha=(\alpha_1,\ldots,\alpha_j), \alpha_i \in \{1,\dots,d\}, j
  \leq k\big\}\,. 
\end{align*}
It is notationally convenient to define the Markov semigroup
$\mathcal{P}_t\colon\C^k\rightarrow \C^k$ associated with
\eqref{eq:main} by
\begin{align}
  (\mathcal{P}_t f)(x) \eqdef \E_x f(X(t))
  \label{eq:markovSemiGP}
\end{align}
where $X(0)=x$ and Markov semigroup $P_h\colon\C^k\rightarrow \C^k$
associated with a single full step of size $h$ of the numerical method
by
\begin{equation*}
  (P_hf)(y) \eqdef\E_y f(Y_1),
\end{equation*}
where $Y_0=y$. Clearly $\|\mathcal{P}_hf\|_0 \leq \|f\|_0$ and
$\|P_hf\|_0 \leq \|f\|_0$. It is also a standard fact, which we
summarize in Appendix~\ref{opBound}, that in our
setting for any $t>0$ and $k \in \N$ if $b, \sigma_1,\dots \sigma_ M
\in \C^k$ then there exists a $C = C(T,k,b,\sigma)$ so that
$\|\mathcal{P}_t f\|_k \leq C \|f\|_k$ is true for all $t \le T$.
 All of these can be rewritten succinctly in the induced operator
norm from $\C^k \rightarrow \C^k$ as $\|\mathcal{P}_t\|_{k\rightarrow
  k} \leq C$, $\|\mathcal{P}_t\|_{0\rightarrow 0} \leq 1$ and
$\|{P}_h\|_{0\rightarrow 0} \leq 1$. Analogously, for any linear
operator $L \colon \C^k \rightarrow \C^\ell$ we will denote the
induced operator norm from $\C^k \rightarrow \C^\ell$ by
$\|L\|_{k\rightarrow \ell}$ which is defined by
\begin{align*}
  \|L\|_{k\rightarrow \ell}= \sup_{f \in \C^k, f\ne 0} \frac{\|L
    f\|_\ell}{\|f\|_k}.
\end{align*}

The following two theorems are the principle results of this article.
They give respectively the weak local and global error of the \algNam
method. 
\begin{theorem}[One-step approximation] \label{thm:local} Assume that
  $b \in \C^6$ and for all $k$, $\sigma_k \in \C^6$ with $\inf_x
  \sigma_k(x) >0$. Then there exists a constant $K$ so that
  \begin{equation}
    \|\mathcal{P}_h - P_h\|_{6\rightarrow 0} \leq K h^3
    \label{eq:lclbnd}
  \end{equation}
  for all $h$ sufficiently small.
\end{theorem}

From this one-step error bound, it is relatively straight-forward to
obtain a global error bound. The following result shows that our
approximation scheme gives a weak approximation of second order.

\begin{theorem}[Global approximation]
  Assume that $b \in \C^6$ and for all $k$, $\sigma_k \in \C^6$ with
  $\inf_x \sigma_k(x) >0$. Then for any $T >0$ there exists a constant
  $C(T)$ such that
  \begin{equation}
    \sup_{0\leq n \le T/h}\|\mathcal{P}_{nh} - P_h^n \|_{6 \rightarrow 0}
    \le C(T) h^2.
    \label{eq:globalbnd} 
  \end{equation}
\end{theorem}

\begin{proof}
  We begin by observing that
  \begin{align*}
    P_h^n- \mathcal{P}_{nh} = \sum_{k=1}^n P_h^{k-1}(P_{h}-
    \mathcal{P}_{h} )\mathcal{P}_{h(n-k)}
  \end{align*}
  and hence since $\sup_{0\leq s \leq
    T}\|\mathcal{P}_{s}\|_{6\rightarrow6} \leq \tilde C(T)$ and
  $\|P_h^k\|_{0\rightarrow 0} \leq 1$, using \eqref{eq:lclbnd} we have
  that for any $n$ with $0 \leq n \leq T/h$
  \begin{align*}
    \|\mathcal{P}_{nh} -P_h^{n}\|_{6 \rightarrow 0}&\leq \sum_{k=1}^n
    \|P_h^{k-1}\|_{0\rightarrow 0}\| P_{h}-
    \mathcal{P}_{h}\|_{6\rightarrow 0}\|\mathcal{P}_{h(n-k)}\|_{6
      \rightarrow 6}\\&\leq \sum_{k=1}^n \tilde C(T)K h^3 =K T\tilde
    C(T) h^2=C(T)h^2.
  \end{align*}
\end{proof}
\begin{remark}
  The restriction that $\inf_x \sigma_k(x) >0$ can likely be relaxed
  if one has some control of the behavior of the solution around the
  degeneracies of $\sigma_k(x)$. This assumption is made to keep the
  proof simple with easily stated assumptions.
\end{remark}
\section{Why the method works}
\label{sec:why}
We now give two different, but related, explanations as to why the \algNamFull
Algorithm is second order accurate in the weak sense.

\subsection{A first point of view}
\label{sec:1stview}
Inserting the expression for $y^*$ from Step 1 of the \algNamFull
Algorithm into Step 2 and disregarding the diffusion terms yields
\begin{align}
  Y_{i} &= Y_{i-1} + h\left[\frac{1}{2\theta}b(y^*) + \left( 1 -
      \frac{1}{2\theta}\right)b(Y_{i-1}) \right] + \dots.
  \label{eq:thewhy1}\\
  &= Y_{i-1} + b(Y_{i-1})h + \frac{b(y^*) - b(Y_{i-1})}{\theta
    h}\frac{h^2}{2} + \dots.
  \label{eq:thewhy2}
\end{align}
Considering \eqref{eq:thewhy1}, we see that when $\theta \approx 1$ we
recover the standard theta method (not to be confused with our use of
$\theta$) with theta $= 1/2$, which is known to be a second order
method for deterministic systems. When $\theta = 1/2$, we recover the
standard trapezoidal or midpoint method.  For $\theta \ne 1/2$, we
simply have a trapezoidal rule where a fractional point of the
interval is used in the construction of the trapezoid.  We will argue
heuristically that the \algNam Algorithm handles the diffusion terms
similarly.  We also note that \eqref{eq:thewhy2} shows that our
algorithm can be understood as an approximation to the two-step Taylor
series where $\theta$ is a parameter used to approximate the second
derivative.  This idea will be revisited in the proof of
Theorem~\ref{thm:local}.

Equation \eqref{eq:main} is distributionally equivalent to
\begin{equation}
  X(t) = X(0) + \int_0^t b(X(s))ds + \sum_{k=1}^M \nu_k \int_0^{\infty}
  \int_0^t 1_{[0,\sigma^2_k(X(s)))}(u) Y_k(du \times ds),
  \label{eq:spacetime}
\end{equation}
where the $Y_k$ are independent space-time white noise
processes\footnote{More precisely, the $Y_k$ are random measures on
  $[0,\infty)^2$ such that if $A,B \subset [0,\infty)^2$ with $A \cap
  B =\emptyset$ then $Y_k(A)$ and $Y_k(B)$ are each independent, mean
  zero Gaussian random variables with variances $\text{Area}(A)$ and
  $\text{Area}(B)$, respectively. Integration with respect to this
  field can be defined in the standard way beginning with adapted
  simple functions which are fixed random variables on fixed
  rectangular sets and then extending by linearity after the
  appropriate It\^o isometry is established.} and all other notation
is as before, in that solutions to \eqref{eq:spacetime} are Markov
processes that solve the same martingale problem as solutions to
\eqref{eq:main}; that is, they have the same generator
(\cite{Kurtz86}).  In order to approximate the diffusion term in
\eqref{eq:spacetime} over the interval $[0,h)$, we must approximate
$Y_k(A_{[0,h)}(\sigma_k^2))$ where $A_{[0,h)}(\sigma_k^2)$ is the
region under the curve $\sigma_k^2(X(t))$ for $0 \le t \le h$.
\begin{figure}
  \subfloat[First step]
  {\label{fig:a}\includegraphics[scale=0.25]{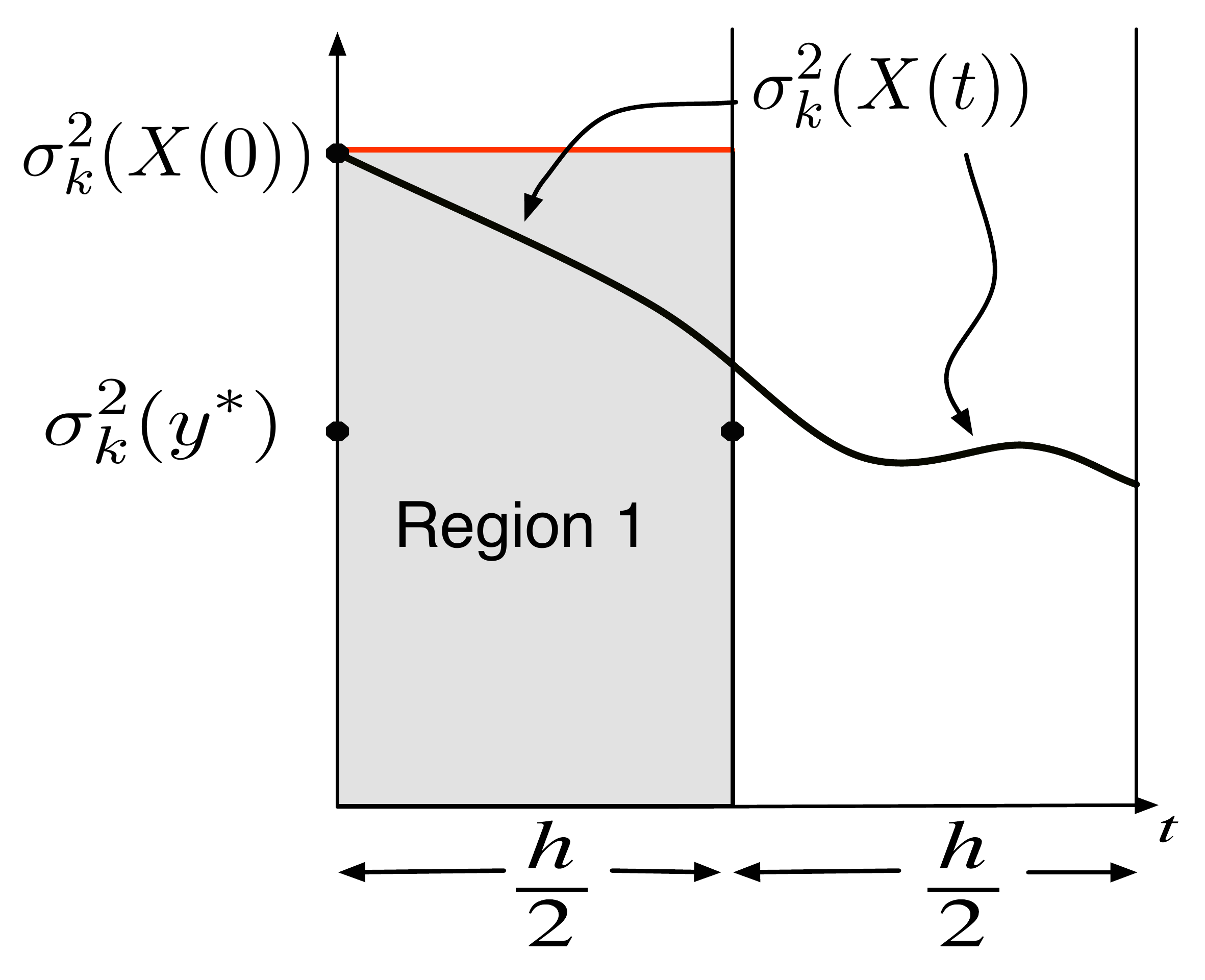}} 
  \subfloat[Desired second step]
  {\label{fig:b}\includegraphics[scale=0.25]{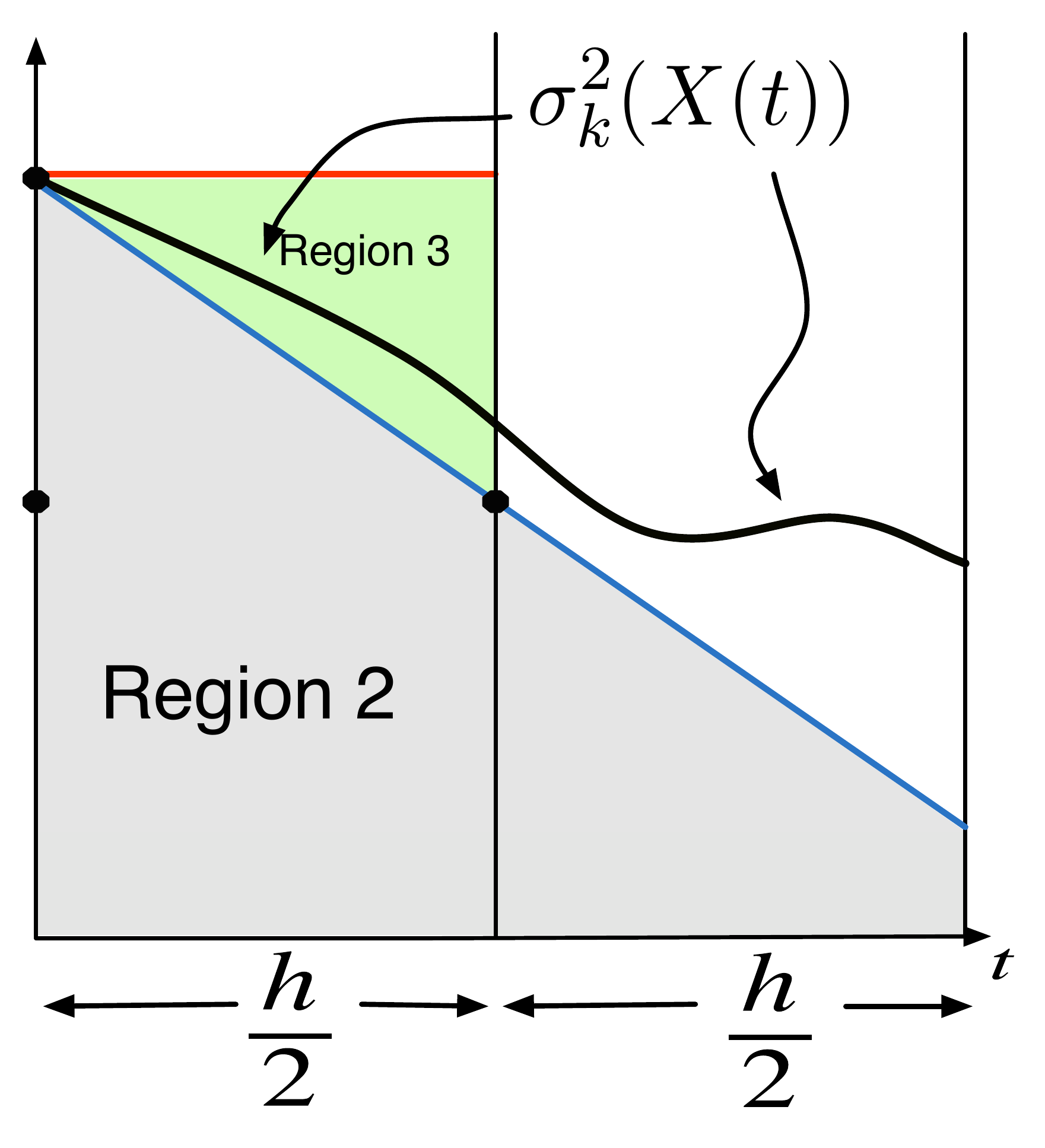}} 
  \subfloat[Used second step]
  {\label{fig:c}\includegraphics[scale=0.25]{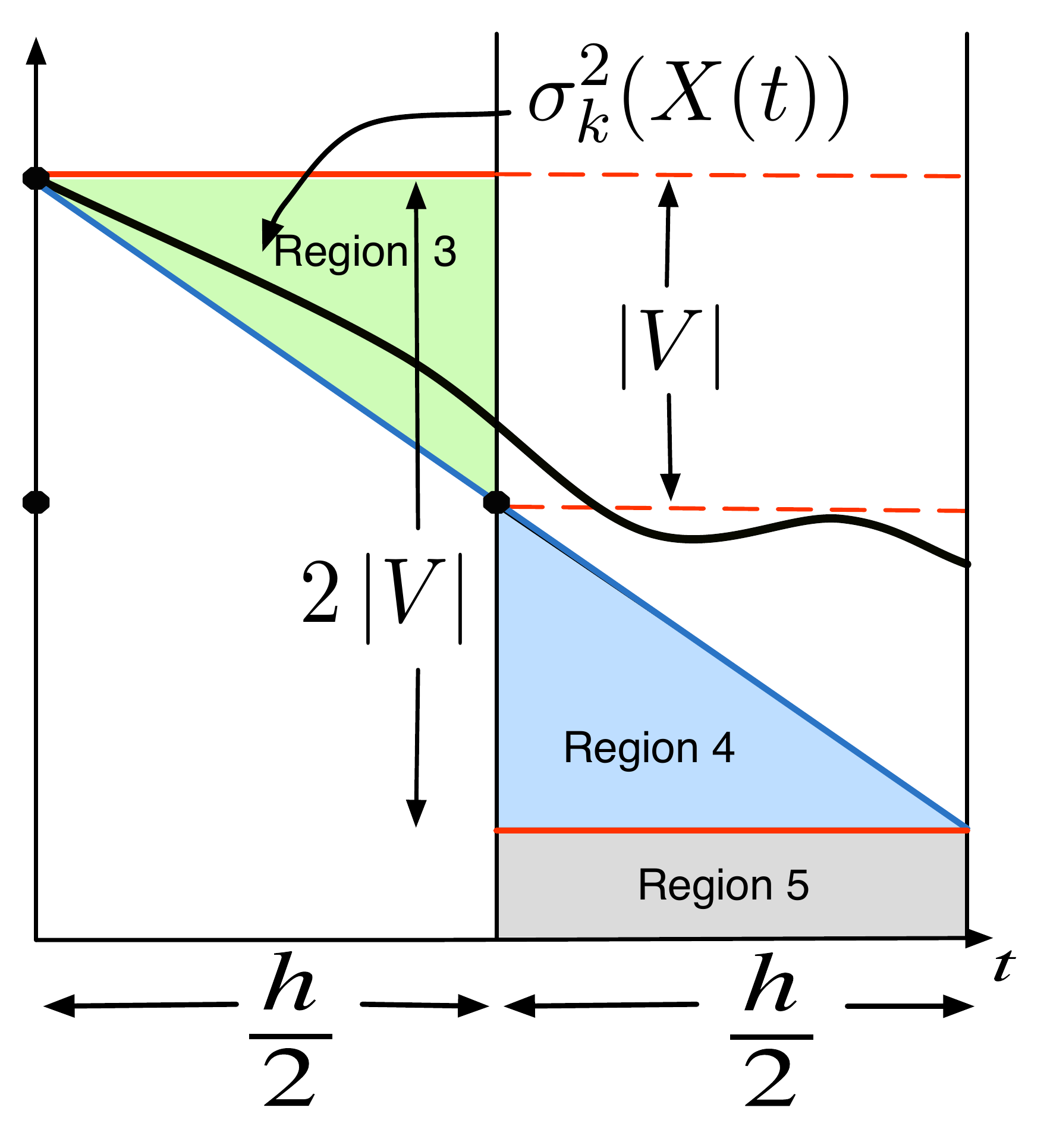}} %
  \caption{A graphical depiction of the \algNam Algorithm with $\theta
    = 1/2$.  In $(a)$ the region of space-time used in the first step
    of the \algNam Algorithm is depicted by the grey shaded Region 1.
    In $(b)$ the desired region to use, in order to perform a
    trapezoidal approximation, would be Region 2. However we have used
    Region 3 in our previous calculation and this is analytically
    problematic to undo.  In $(c)$, where $V = \sigma^2_k(y^*) -
    \sigma^2_k(X(0))$, we see that Region 5 gives the correct amount
    of new area wanted as subtracting off the area of Region 4
    ``offsets'' the used area of Region 3.  The case $\theta \ne 1/2$
    is similar.}
  \label{fig:exp}
\end{figure}

We consider a natural way to approximate $X(h)$ and focus on the
double integral in \eqref{eq:spacetime} for a single $k$.  We also
take $\theta = 1/2$ for simplicity and simply note that the case
$\theta \ne 1/2$ follows similarly.  We begin by approximating the
value $X(h/2)$ by $y^*$ obtained via an Euler approximation of the
system on the interval $[0, h/2)$.  To do so, we hold $X(t)$ fixed at
$X(0)$ and see that we need to calculate $Y_k(\text{Region 1})$, where
Region 1 is the grey shaded region in Figure \ref{fig:exp}(a).
Because
\begin{equation*}
  Y_k(\text{Region 1}) \ \eqdist \ N(0,\sigma_k^2(X(0))  h/2) \
  \eqdist \ \sigma_k(X(0))\sqrt{\tfrac{h}{2}} \ N(0,1), 
\end{equation*}
we see that this step is equivalent in distribution to Step 1 of
Algorithm \ref{eq:main}. (Here and in the sequel, `` $\eqdist$ '' denotes
``equal in distribution.'')

If we were trying to determine the area under the curve
$\sigma_k^2(X(t))$ using an estimated midpoint $y^*$ for a
deterministic $X(t)$, one natural (and common) way would be to use the
area of Region 2, where Region 2 is the grey shaded region in Figure
\ref{fig:exp}(b).  Such a method would be equivalent to the
trapezoidal rule given in \eqref{eq:thewhy1}.  However, in our setting
we would have to ignore, or subtract off, the area already accounted
for in Region 3, which is depicted as the shaded green section of
Figure \ref{fig:exp}(b).  In doing so, the random variable needed in
order to perform this step would necessarily be dependent upon the
past (via Region 3), and our current analysis would break down.
However, noting that Region 3 has the same area as Region 4, as
depicted by the blue shaded region in Figure \ref{fig:exp}(c), we see
that it would be reasonable to expect that if one only uses Region 5,
as depicted as the grey shaded region in Figure \ref{fig:exp}(c), then
the accuracy of the method should be improved as we have performed a
trapezoidal type approximation.  Because
\begin{align*}
  Y_k(\text{Region 5}) \ \eqdist \ N\Big( 0 , \left(\sigma_k^2(X(0))
      + 2V\right) \tfrac{h}{2}\Big) &\eqdist \,\sqrt{\sigma_k^2(X(0))
    + 2V}\sqrt{\tfrac{h}{2}} N(0,1)\\ &= \sqrt{2\sigma_k^2(y^*) -
    \sigma_k^2(X(0))}\sqrt{\tfrac{h}{2}} N(0,1),
\end{align*}
where $V = \sigma_k^2(y^*) - \sigma_k^2(X(0))$, we see that this is
precisely what is carried out by Step 2 of the \algNamFull
Algorithm.

\subsection{A second point of view}
To obtain a higher order method one must both approximate well the
expected drift term as well as the quadratic variation of the process.
The basic idea of the \algNam Algorithm is to make a preliminary step
using an Euler approximation and then use this step to make a higher
order approximation to the drift integral and to the quadratic
variation integral.  Similar to \eqref{eq:thewhy1} the desired one
step approximation to the quadratic variation integrals are
\begin{align*}
  \int_0^h \sigma_k^2(X(s)) ds \approx
  h\left[\frac{1}{2\theta}\sigma_k^2(y^*) + \left( 1 -
      \frac{1}{2\theta}\right)\sigma_k^2(Y_{i-1}) \right],
\end{align*}
where all notation is as before.

Considering just the variance terms of the quadratic variation, we let
$\{e_i\}$ be an orthonormal basis and see that our method yields the
approximation
\begin{align*}
  \text{Var}(X(h)\cdot e_i )&\approx
  \sum_{k=1}^M\text{Var}\Big(\sigma_k( Y_{0} )(\nu_k\cdot
    e_i)\eta_{1k} \sqrt{\theta h} + \sqrt{\big[\alpha_1
      \sigma^2_k(\ym) - \alpha_2 \sigma^2_k(Y_{0})\big]^+} (\nu_k\cdot
    e_i) \eta_{2k} \sqrt{(1-\theta)h}\Big)\\&=
    \sum_{k=1}^M\E\Big(\sigma_k^2( Y_{0} )\theta 
  +\big[\alpha_1\sigma^2_k(\ym) - \alpha_2\sigma^2_k(Y_{0})\big]^+ (1
  - \theta)\Big) (\nu_k\cdot e_i)^2h.
\end{align*}
If the step-size is sufficiently small then, $\big[\alpha_1
\sigma^2_k(\ym) - \alpha_2 \sigma^2_k(Y_{0})\big]^+$ is positive with
high probability because of our uniform ellipticity assumption; and
hence,
\begin{align*}
  \text{Var}(X(h)\cdot e_i )\approx \E \sum_{k=1}^M (\nu_k\cdot
  e_i)^2\; \Big(\frac{1}{2\theta}\sigma_k^2(y^*) + \Big( 1 -
      \frac{1}{2\theta}\Big)\sigma_k^2(Y_{i-1})
  \Big) h
\end{align*}
which is a locally third order approximation to the true quadratic
variation integral of
\begin{equation*}
  \text{Var}(X(h)\cdot e_i) = \E \sum_{k=1}^M (\nu_k\cdot e_i)^2\int_0^h
  \sigma_k^2(X(s)) ds. 
\end{equation*}
Notice that it was important in this simple analysis that the
direction of variation $\nu_k$ stayed constant over the interval so
that the two terms could combine exactly.
Of course, one should really be computing the full quadratic
variation, including terms such as $\text{Cov}(X(h)\cdot e_i,
X(h)\cdot e_j)$, but they follow the same pattern as above because
each is a linear combination of the integral terms $\int_0^h
\sigma_k^2(X(s))ds$.
\section{Proof of Local Error Estimate}
\label{sec:proof}
We now give the proof of the local error estimate given in Theorem
\ref{thm:local} which is the central result of this paper.
\begin{proof}(of Theorem~\ref{thm:local}) We need to show that there
  exists a constant $K$ so that for any $f \in \C^6$ one has
  \begin{align*}
    \left |\E f(Y_1) - \E f(X(h)) \right| \leq K \|f\|_6 h^3\, .
  \end{align*}
  Hence for the reminder of the proof we fix an arbitrary $f \in
  \C^6$.
  Observe that Step 1 of the \algNam Algorithm produces a value,
  $\ym$, that is distributionally equivalent to $y(\theta h)$, where
  $y(t)$ solves
  \begin{align}
    dy(t) = b(y(0))dt + \sum_{k = 1}^M \sigma_k(y(0)) \;\nu_k \
    dW_k(t), \quad y(0) = x_0.
    \label{eq:part1}
  \end{align}
  Likewise, Step 2 of the \algNam Algorithm produces a value, $Y_1$,
  that is distributionally equivalent to $y(h)$, where $y(t)$ solves
  \begin{align}
    dy(t) = (\alpha_1 b(\ym) - \alpha_2 b(x_{0}))dt + \sum_{k=1}^M
    \sqrt{[\alpha_1 \sigma^2_k(\ym) - \alpha_2\sigma^2_k(x_0)]^+}
    \;\nu_k \ dW_k(t), \quad y(\theta h) = \ym.
    \label{eq:part2}
  \end{align}
  Let $\Ft$ denote the filtration generated by the Weiner processes
  $W_k(t)$ in \eqref{eq:part1} and \eqref{eq:part2}.  Then,
  \begin{equation}
    \E f(y(h)) = \E\,[\,\E [ f(y(h))\,|\,\Ftt ]\,] \eqdef
    \E\,[\,\E_{\theta h} f(y(h))],  
    \label{eq:condition}
  \end{equation} 
  where we have made the definition $\E_{\theta h}[\ccdot]\eqdef
  \E[ \ \cdot \ | \ \Ftt]$.

  Let $A$ denote the generator for the process \eqref{eq:main}, $B_1$
  denote the generator for the process \eqref{eq:part1}, and $B_2$
  denote the generator for the process \eqref{eq:part2} conditioned
  upon $\Ftt$.  Then
  \begin{align*}
    (Af)(x) &= f'[b](x) + \frac{1}{2} \sum_k \sigma_k^2 f''[\nu_k,
    \nu_k](x)\\
    (B_1f)(x) &= f'[b(x_0)](x) + \frac{1}{2} \sum_k \sigma_k(x_0)^2
    f''[\nu_k, \nu_k](x)\\
    (B_2 f)(x) &= f'[\alpha_1 b(\ym) - \alpha_2 b(x_0)](x) +
    \frac{1}{2} \sum_k [\alpha_1 \sigma_k(\ym)^2 - \alpha_2
    \sigma_k(x_0)^2]^+ f''[\nu_k, \nu_k](x),
  \end{align*}
  where $f'[\xi](z)$ is the derivative of $f$ in the direction $\xi$
  evaluated at the point $z$.  Note that $(Af)(x_0) = (B_1f)(x_0)$.
  For any integer $k \ge 2$ we define recursively $(A^k f)(x)\eqdef
  (A(A^{k-1}f))(x),$ and similarly for $B_1$ and $B_2$.  By repeated
  application of the It\^o-Dynkin formula, see \cite{Oksendal03}, we
  have
  \begin{align}
    \E_{\theta h} f(y(h)) &= f(\ym) + \int_{\theta h}^h \E_{\theta h}
    (B_2 f)(y(s)) \ ds \notag\\
    &= f(\ym) + (B_2f)(\ym)(1 - \theta)h + \int_{\theta h}^h
    \int_{\theta h}^s
    \E_{\theta h} (B_2^2 f)(y(r))\ dr \ ds \notag \\
    \begin{split}  
      &= f(\ym) + (B_2f)(\ym)(1 - \theta)h + (B_2^2 f)(\ym)
      \frac{(1 - \theta)^2 h^2}{2}\\
      & \hspace{.6in} + \int_{\theta h}^h \int_{\theta h}^s
      \int_{\theta h}^r \E_{\theta h} (B_2^3 f)(y(u))\ du \ dr \ ds.
    \end{split}
    \label{eq:1stexpansion}
  \end{align}
  The term $(B_2^3 f)(y(u))$ depends on the first six derivatives of
  $f$.  Therefore, since $f \in \C^6$
  \begin{equation}
    \left| \int_{\theta h}^h \int_{\theta h}^s \int_{\theta h}^r
      \E_{\theta h} (B_2^3 f)(y(u))\ du \ dr \ ds \right| \le C \|f\|_6 h^3,
    \label{eq:trivial_bnd}  
  \end{equation}
  for some constant $C$.  Combining \eqref{eq:condition},
  \eqref{eq:1stexpansion}, \eqref{eq:trivial_bnd}, and recalling that
  $\E f(Y_1)= \E f(y(h))$ gives
  \begin{align}
    \E f(Y_1)&= \E \ [ \ \E_{\theta h} f(y(h))] = \E \ f(\ym) + \E \
    (B_2f)(\ym)(1 - \theta)h + \E \ (B_2^2 f)(\ym)\frac{(1 - \theta)^2
      h^2}{2} + O(h^3).
    \label{eq:partway}
  \end{align}
  Here and in the sequel, we will write $F = G + O(h^p)$ to mean that
  there exist a constant $K$ depending on only $\sigma$ and $b$ so
  that for all initial conditions $x_0$
  \begin{align}\label{eq:Onotation}
    |F - G| \leq K \|f\|_6 h^p\, ,
  \end{align}
  for $h$ sufficiently small. In the spirit of the preceding
  calculation, repeated application of the It\^o-Dynkin formula to
  \eqref{eq:main} produces
  \begin{equation*}
    \E f(X(h)) =   f(x_0) + (Af)(x_0)h + (A^2f)(x_0)\frac{h^2}{2} + O(h^3).
  \end{equation*}  
  The proof of the theorem is then completed by
  Lemma~\ref{lem:keyEstimate} given below. Its proof, which is
  straightforward but tedious, is given in the appendix.
\end{proof}

\begin{lemma} \label{lem:keyEstimate}%
  Under the assumptions of Theorem~\ref{thm:local}, for all $h>0$
  sufficiently small and $f \in \C^6$ one has
  \begin{align*}
    \E \left[ f(\ym) + (B_2f)(\ym)(1 - \theta)h + (B_2^2 f)(\ym)
      \frac{(1 - \theta)^2 h^2}{2} \right] =& f(x_0) + (Af)(x_0)+
    (A^2f)(x_0)\frac{h^2}2 + O(h^3)\,.
  \end{align*}
\end{lemma}

\begin{remark}
  Comparing equation \eqref{eq:thewhy2} and Lemma
  \ref{lem:keyEstimate} shows that our algorithm can be viewed as
  providing an approximation to the two step Taylor series
  approximation.
\end{remark}

\section{Examples}
\label{sec:examples}
We present two examples that demonstrate the rate of convergence of
the \algNam Algorithm with $\theta = 1/2$.  In each example we shall
compare the accuracy of the proposed algorithm to that of Euler's
method and a ``midpoint drift'' algorithm defined via repetition of
the following steps
\begin{align}
  \begin{split}
    \ym &= Y_{i-1} + b(Y_{i-1}) \frac{h}{2}\\
    Y_i &= Y_{i-1} + b(\ym)h + \sum_{k=1}^M \sigma_k(Y_{i-1})\nu_k\;
    \eta_k\; \sqrt{h},
  \end{split}
  \label{eq:mdptdrift}
\end{align} 
where the notation is as before.  We compare the proposed algorithm to
that given via \eqref{eq:mdptdrift} to point out that the gain in
efficiency being demonstrated is not solely due to the fact that we
are getting better approximations to the drift terms,
but also because of the superior approximation of the diffusion terms.

\subsection{First Example.}
 \label{ex:OU}
\begin{figure}
    \begin{center}
      \subfloat[First Example] {\label{fig:OU}
        \includegraphics[height=2.25in]{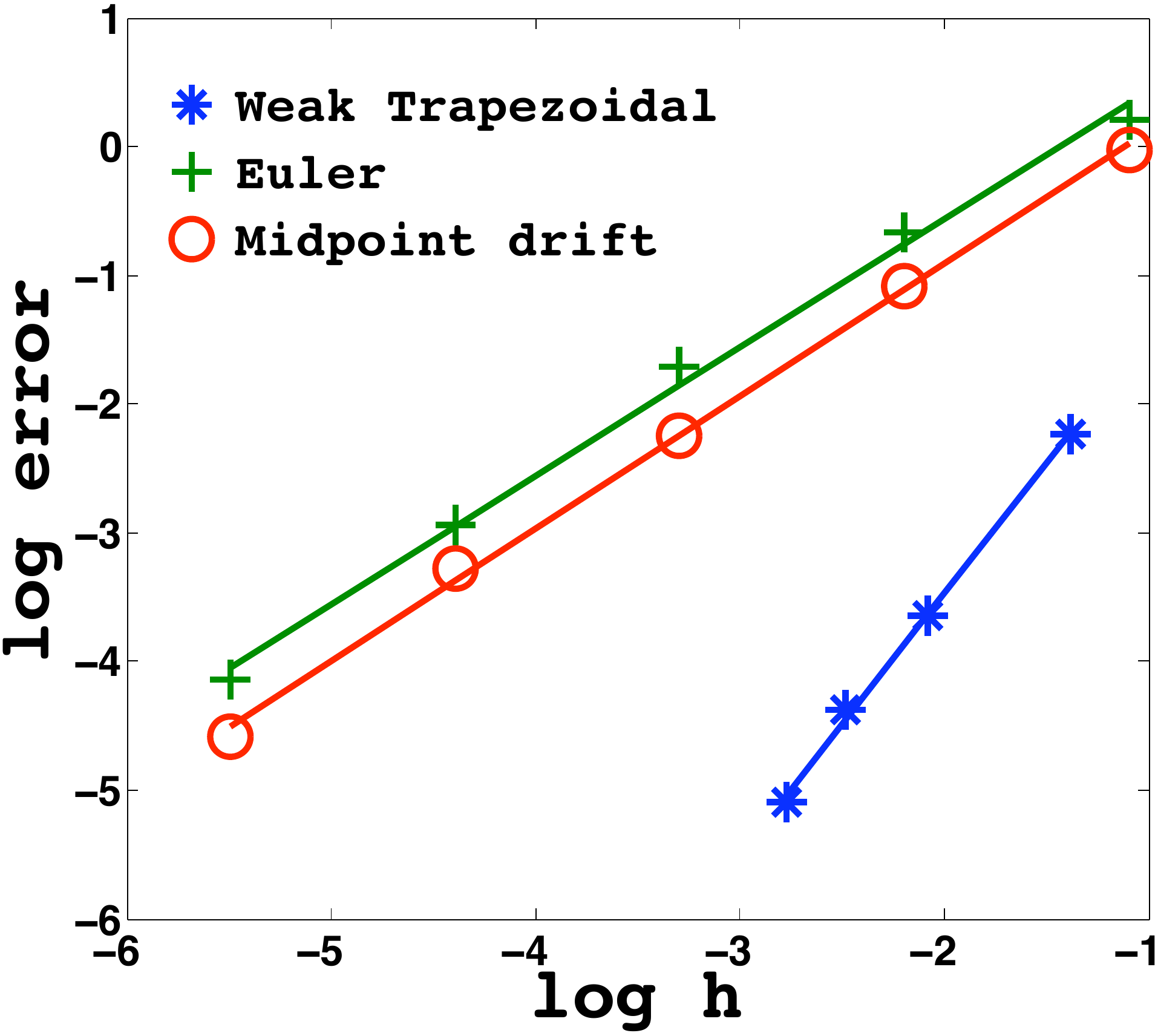}} \qquad\qquad
      \subfloat[Second Example]
      {\label{fig:talay}\includegraphics[height=2.25in]{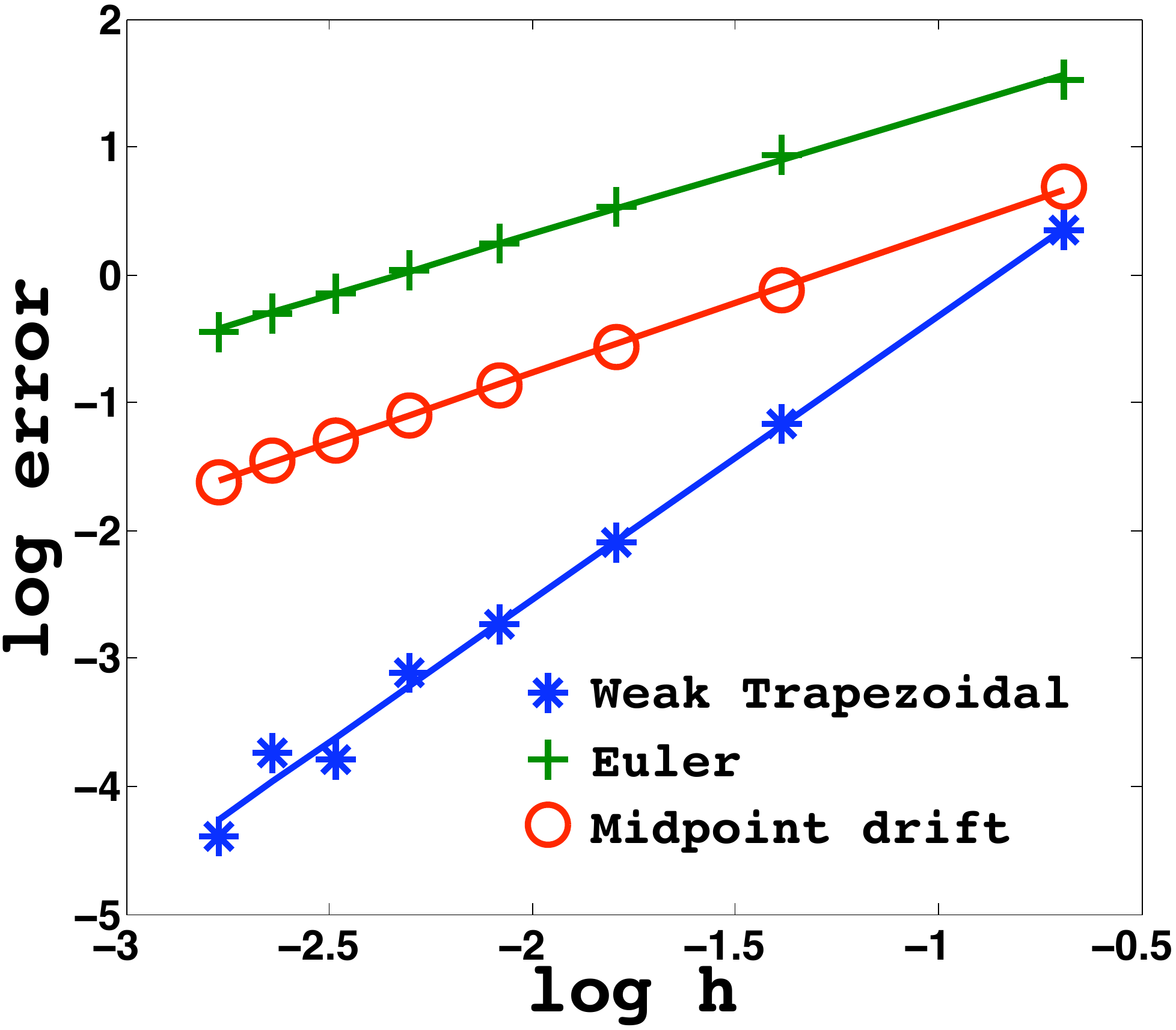}}
    \end{center} %
    \caption{ Log-log plots of the step-size versus the error for the
      three different algorithms.  In (a) the example \eqref{eq:OU} is
      considered.  The best fit lines for the data (shown) have slopes
      2.029, .998, and 1.030, for the \algNam Algorithm, Euler's
      method, and the midpoint drift method, respectively.  In (b) the
      example in \eqref{eq:talay} is considered.  The best fit lines
      for the data (shown) have slopes 2.223, .952, and 1.098, for
      the \algNam Algorithm, Euler's method, and the midpoint drift
      method, respectively.  In both examples all results agree with
      what was expected.}
 \end{figure}

  Consider the system 
  \begin{equation}
    \left[ \begin{array}{c}
        dX_1(t) \\
        dX_2(t)
      \end{array}\right]
    = \left[ \begin{array}{c}
        X_1(t) \\
        0
      \end{array}\right] + X_1(t) \left[ \begin{array}{c}
        0 \\
        1
      \end{array}\right]dW_1(t) + \frac{1}{10}\left[ \begin{array}{c}
        1 \\
        1
      \end{array}\right]dW_2(t),
    \label{eq:OU}
  \end{equation}
  where $W_1(t)$ and $W_2(t)$ are standard Weiner processes.  In our
  notation $b_1(x) = x_1$, $b_2(x) = 0$, $\sigma_1(x) = x_1$,
  $\sigma_2(x) = 1/10$, and $\nu_1 = [0,1]^T$, $\nu_2 = [1,1]^T$.
  Note that the noise does not commute.  It is an exercise to show
  that
  \begin{equation} 
    \E X_2(t)^2 = \E \ X_2(0)^2 - \frac{1}{2}\E \ X_1(0)^2 +
    \frac{1}{400}e^{2t}(200 \E X_1(0)^2 + 1) + 
    \frac{t}{200}  - \frac{1}{400}.
    \label{eq:OUsol}
  \end{equation}
  For both Euler's method and the midpoint drift method
  \eqref{eq:mdptdrift} we used step sizes $h_k = 1/3^k$, $k \in
  \{1,2,3,4,5\}$ and initial condition $X_1(0) = X_2(0) = 1$ to
  generate $500,000$ sample paths of the system \eqref{eq:OU}.  We
  then computed
  \begin{equation}\label{eq:error}
    \text{error}_k(t) = \E X_2(t)^2 - \frac{1}{5\times 10^5} \sum_{i =
      1}^{ 5\times 10^5} \overline X_2^{h_k}(t)^2,
  \end{equation}
  where $\overline X^{h_k}(t)$ is the sample path generated
  numerically and $\E X_2(t)^2$ is given via \eqref{eq:OUsol}.  We
  also generated $10,000,000$ sample paths using the \algNam Algorithm
  with the same initial condition and step sizes $h_k = 1/(4k)$, $k
  \in \{1,2,3,4\}$.  We then computed $\text{error}_k(t)$ similarly to
  before.  The outcome of the numerical experiment is summarized in
  Figure \ref{fig:OU} where we have plotted $\log(h_k)$ versus
  $\log(|\text{error}_k(1)|)$ for the different algorithms.  As
  expected, we see that the \algNam Algorithm gives an error that
  decreases proportional to $h^2$, whereas the other two algorithms
  give errors that decrease proportional to $h$.
    
\subsection{Second Example.}
%
Now consider the following system that is similar to one considered in
\cite{talay1990}
\begin{align}
  \begin{split}
    \left[ \begin{array}{c}
        dX_1(t) \\
        dX_2(t)
      \end{array}\right]
    = \left[ \begin{array}{c}
        -X_2(t) \\
        X_1(t)
      \end{array}\right] &+ \sqrt{\frac{\sin^2(X_1(t) + X_2(t)) +
        6}{t + 1}} \left[ \begin{array}{c}
        1 \\
        0
      \end{array}\right]dW_1(t)\\
    & + \sqrt{\frac{\cos^2(X_1(t) + X_2(t)) + 6}{t + 1}}\left[
      \begin{array}{c}
        0 \\
        1
      \end{array}\right]dW_2(t),
  \end{split}
  \label{eq:talay}  
\end{align}
where $W_1(t)$ and $W_2(t)$ are independent Weiner processes.  It is
simple to show that
\begin{equation}
  \E |X(t)|^2 = \E X(0)^2 + 13\log(1 + t).
  \label{eq:talaysol}
\end{equation}  We used step sizes $h_k = 1/(2k)$, $k \in \{1,2,\dots,8\}$, to
generate five million approximate sample paths of the system
\eqref{eq:talay} using each of: (a)  \algNam Algorithm, (b)
Euler's method, and (c) the midpoint drift method
\eqref{eq:mdptdrift}.  We then computed
\begin{equation*}
  \text{error}_k(t) = \E |X(t)|^2 - \frac{1}{5\times 10^6} \sum_{i =
    1}^{5\times 10^6} |\overline X^{h_k}(t)|^2,
\end{equation*}
where $\overline X^{h_k}(t)$ is the sample path generated numerically
and $\E |X(t)|^2$ is given via \eqref{eq:talaysol}.  The outcome is
summarized in Figure \ref{fig:talay} where we have plotted $\log(h_k)$
versus $\log(|\text{error}_k(1)|)$ for the different algorithms.  As
before, we see that the \algNam Algorithm gives an error that
decreases proportional to $h^2$, whereas the other two algorithms give
errors that decrease proportional to $h$.
  
\vspace{.1in}

\begin{remark}
  We note that in both examples we needed to average over an extremely
  large number of computed sample paths in order to estimate
  error$_k(t)$ for the \algNam Algorithm.  This is due to the fact
  that the increased accuracy of the method quickly makes sampling
  error the dominant error.
\end{remark}

\section{The effect of varying $\theta$}
\label{sec:theta}
\begin{figure}
\begin{center}
\subfloat[\% of degenerate steps vs
  $\theta$]{\label{fig:theta}\includegraphics[height=2.25in]{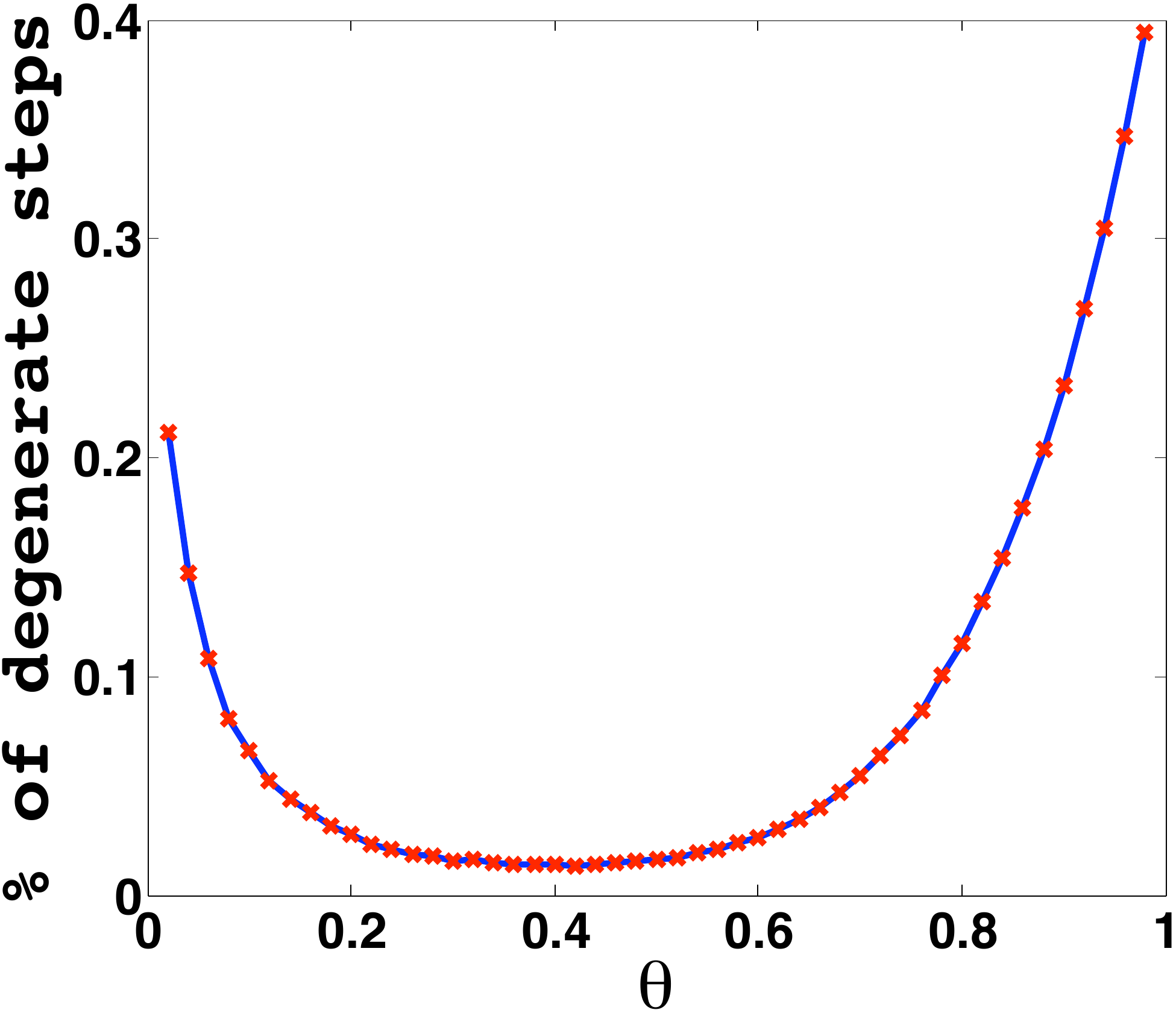}}
\qquad 
\subfloat[accuracy for different
$\theta$]{\label{fig:loglogtheta}\includegraphics[height=2.25in]{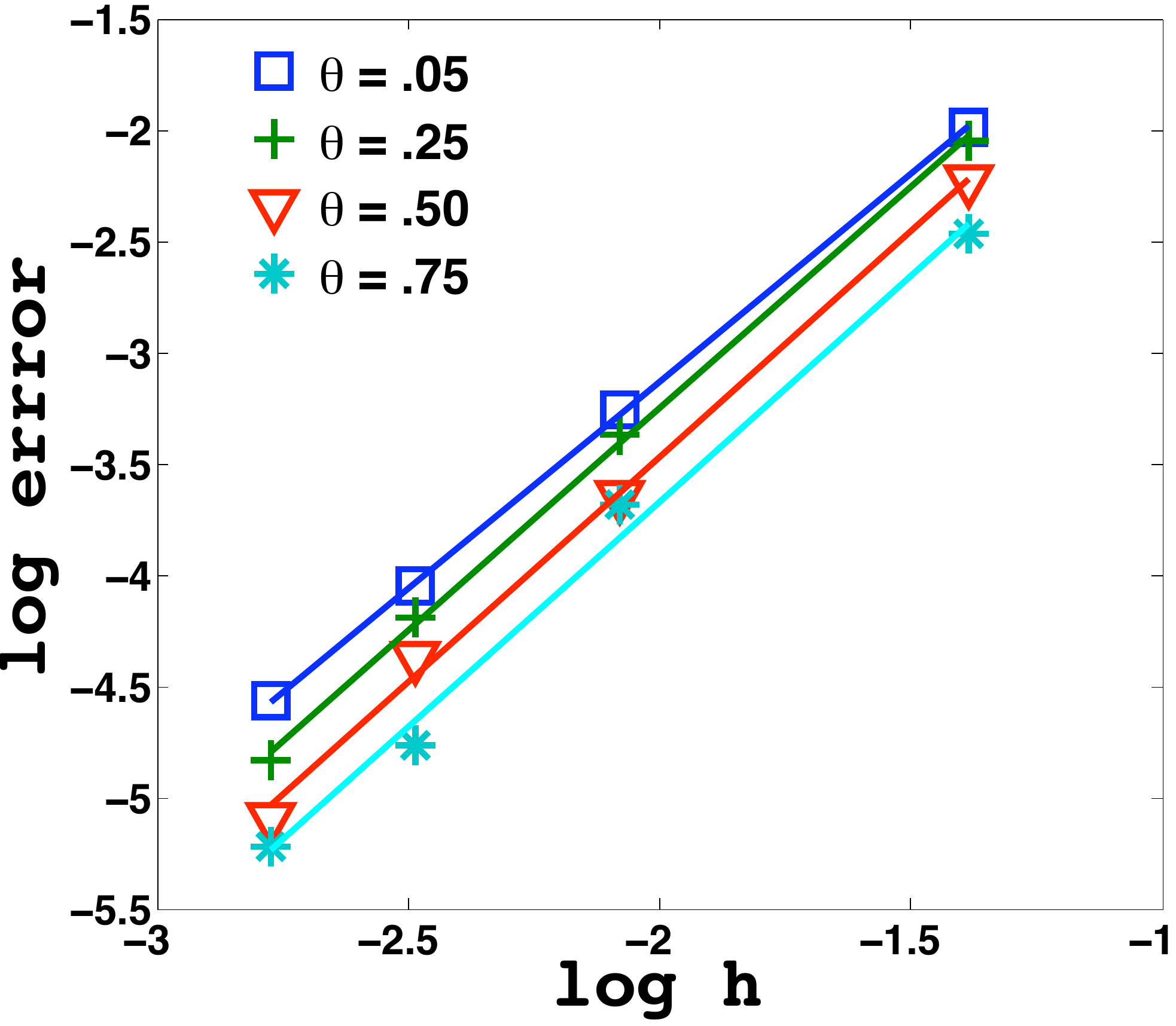}}
\caption{(a) The number of degenerate steps for the \algNam Algorithm
  applied to \eqref{eq:testsys} with $h = 1/10$ and different values
  of $\theta$.  (b) The $\log h$ vs $\log(|\textrm{error}|)$ plot is
  given for different choices of $\theta$ for the \algNam Algorithm
  applied to \eqref{eq:OU} where the \textrm{error} is defined
  similarly to \eqref{eq:error}.  The best fit lines for the data
  (shown) have slopes 1.865, 1.996, 2.029, and 2.033 for $\theta =
  .05, .25, .50, .75$, respectively.}
\end{center}
\end{figure}

The term $\big[\alpha_1 \sigma_k^2(y^*) - \alpha_2
  \sigma_k^2(Y_{i-1})\big]^+$ in Step 2 of the \algNam Algorithm will
yield zero, and the given step will have a local error of only
$O(h^2)$, if 
\begin{equation*}
  \alpha_1 \sigma_k^2(y^*) < \alpha_2 \sigma_k^2(Y_{i-1})
  \Longleftrightarrow  \sigma_k^2(y^*) \, < \, \frac{\alpha_2}{\alpha_1}
  \sigma_k^2(Y_{i-1}) \, = \, (1 - 2\theta +
  2\theta^2)\sigma_k^2(Y_{i-1}). 
\end{equation*}
We will call such a step a ``degenerate'' step.  The function
$g(\theta) = 1 - 2\theta + 2\theta^2$ is minimized at $g(1/2) = 1/2$,
and $g(\theta) \to 1$ as $\theta \to 0$ or $\theta \to 1$.  Therefore,
as mentioned Remark~\ref{rem:theta}, one would expect that as $\theta
\to 0$ or $\theta \to 1$ more steps will be degenerate, and a decrease
in accuracy, together with a bias against $\sigma_k$ decreasing, would
follow.  Using a step-size of $h = 1/10$, we tracked the percentage of
degenerate steps for the simple system
\begin{equation}
  dX(t) = \sqrt{X(t)^2 + 1}\ dW(t), \qquad X(0) = 1,
  \label{eq:testsys}
\end{equation}
where $W(t)$ is a standard Weiner process.  To do so, we computed
$10,000$ sample paths over the time interval $[0,1]$ for each of
$\theta = .02k$, \, $k\in \{1,\dots,49\}$.  The results are shown in Figure
\ref{fig:theta} where the behavior predicted above is seen.  Curiously,
the minimum number of rejections takes place at $\theta = .42$.  
It is also worth noting that one can check on computer software that
in the general case the coefficient of $h^3$ for the one-step error
grows like $1/\theta$ as $\theta \to 0$.  This does not happen in the
deterministic case \eqref{eq:thewhy1}.

While the above considerations give some interesting insight into the
effect of various $\theta$, the situation is more complex.  A $\theta$
closer to one should give the method more stability, albeit at an
expense as the rejection fraction increases as $\theta$ approaches
one.  It would be interesting to perform a stability analysis in the
spirit of \cite{Higham00} to better understand the effect of $\theta$.
In lieu of this, Figure~\ref{fig:loglogtheta} gives the result of a
convergence analysis of the \algNam Algorithm applied to \eqref{eq:OU}
with different choices of $\theta$. Interestingly, larger $\theta$
seem to result in smaller (and hence better) convergence rate
prefactors.  This seems to indicate that in at least this example
stability is an issue. 

The performance of the \algNam Algorithm as a function of $\theta$ is
a topic deserving further consideration, but combining the above shows
that $\theta = 1/2$ is a reasonable first choice, though stability
considerations might lead one to consider a $\theta$ closer to 1.

\section{Comparison to Richardson Extrapolation}\label{sec:richardson}
It is illustrative to compare the \algNam Algorithm to Richardson
extrapolation, which from a certain point of view is the method in the literature that
is most similar to ours.  See \cite{talay1990} for complete details of
Richardson extrapolation in the SDE setting.

Let $Z_{h/2}(t)$ and $Z_h(t)$ denote approximate sample paths of
\eqref{eq:main} generated using Euler's method with step sizes of
$h/2$ and $h$, respectively.  For all $f$ satisfying mild assumptions,
both $\E f(Z_{h/2}(t))$ and $\E f(Z_h(t))$ will approximate $\E
f(X(t))$ with an order of $O(h)$.  However, Richardson extrapolation
may be used and the linear combination $2\E f(Z_{h/2}(t)) - \E f(Z_h)$
will approximate $\E f(X(t))$ with an order of $O(h^2)$ (see
\cite{talay1990} ).  Of course, taking $f$ to be the identity shows that
the linear combination $2Z_{h/2}(t) - Z_h(t)$ gives an $O(h^2)$
approximate of the mean of the process.  As Richardson extrapolation
does not compute a single path, but instead uses the statistics from
two to achieve a higher order of approximation for a given statistic,
we will compare one step of the \algNam Algorithm with a step-size of
$h$, to one step of size $h$ of the process $2Z_{h/2}(t) - Z_h(t)$
with the clear understanding that $2Z_{h/2}(t) - Z_h(t)$ is only
$O(h)$ accurate for higher moments.

\begin{figure}
  \begin{center}
  \subfloat[When the process increases]
  {\label{fig:RE_a}\includegraphics[scale=0.28]{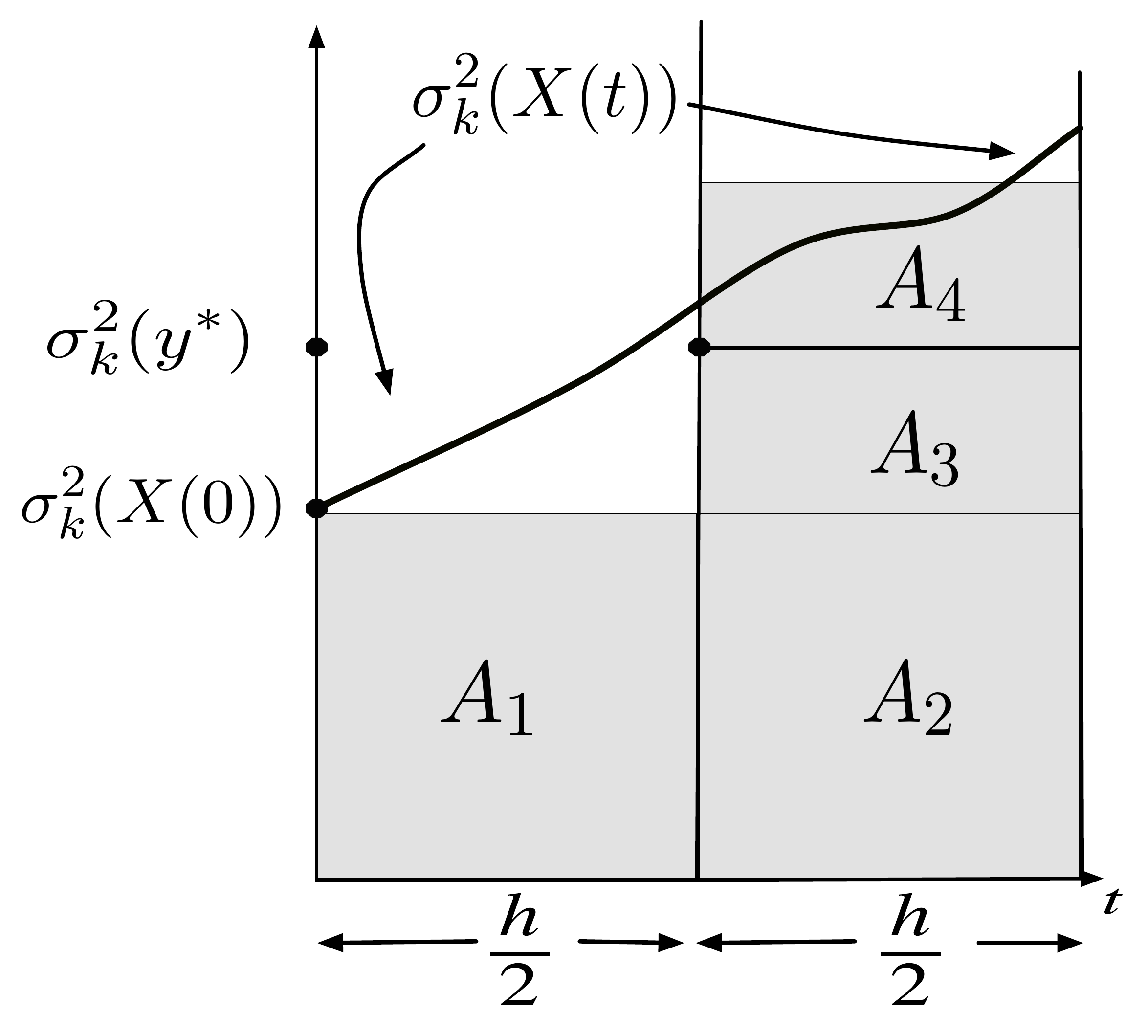}}
  \subfloat[When the process decreases]
  {\label{fig:RE_b}\includegraphics[scale=0.28]{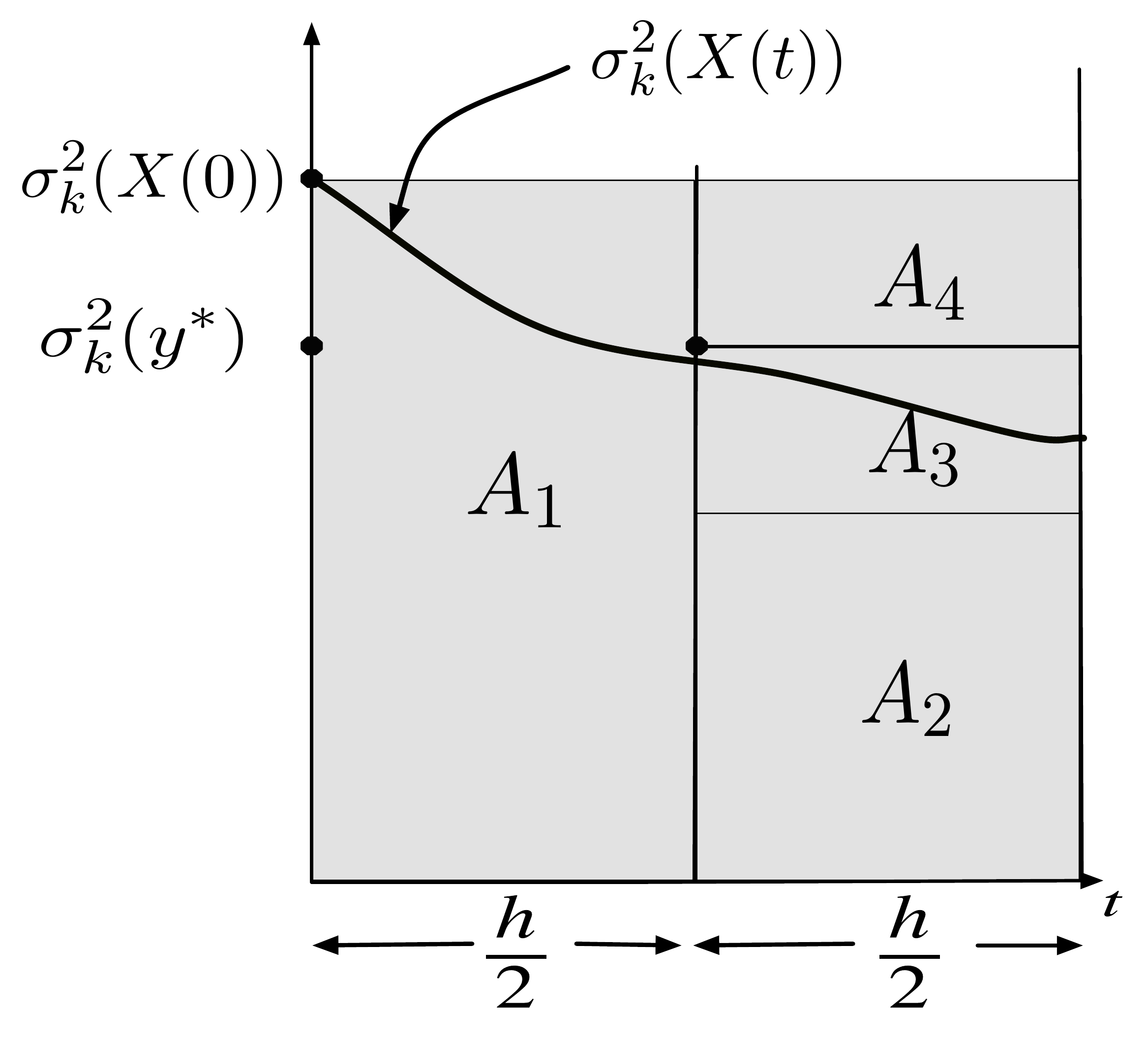}}
  \end{center}
  \caption{The areas of space-time utilized by $2Z_{h/2} - Z_h$ and
    the \algNam Algorithm for a single $k$ and a single step. In
    \ref{fig:RE}(a), $\sigma_k^2(X(t))$ increases and $2Z_{h/2} - Z_h$ uses
    $\eta_{A_1} + \eta_{A_2} + 2\eta_{A_3}$, whereas the \algNam
    Algorithm uses $\eta_{A_1} + \eta_{A_2} + \eta_{A_3} +
    \eta_{A_4}$.  In the case when $\sigma_k^2(X(t))$ decreases,
    \ref{fig:RE}(b) above, the processes use $\eta_{A_1} + \eta_{A_2}
    + \eta_{A_3} - \eta_{A_4}$ and $\eta_{A_1} + \eta_{A_2}$,
    respectively.  In both cases, it is the better use of the areas by
    the \algNam Algorithm that achieves a higher order of
    convergence.}
  \label{fig:RE}
\end{figure}

Recall that systems of the form \eqref{eq:main} are equivalent to
those driven by space-time white noise processes \eqref{eq:spacetime}.
As in Section \ref{sec:1stview}, we consider how each method uses the
areas of $[0,\infty)^2$ associated to $Y_k(du \times ds)$ from
\eqref{eq:spacetime}  during one step. We will proceed considering a
single $k$ since it is sufficient to illustrate the point.  For
$A_i \subset [0,\infty)^2$, we denote by $\eta_{A_i}$ a normal random
variable with mean $0$ and variance area$(A_i)$. Recall that
$\eta_{A_i}$ and   $\eta_{A_j}$  are independent as long  as $A_i \cap
A_j$ has Lebesgue measure zero. Consider
\eqref{fig:RE}(a) in which we are supposing that $\sigma_k^2(X(t))$
increases over a single time-step.  The change in the process
$Z_{h/2}$ due to this $k$ would be $\nu_k$ times
\begin{equation*}
  \eta_{A_1} + \eta_{A_2} + \eta_{A_3}.	
\end{equation*}
Similarly, the change in $Z_h$ would be $\nu_k$ times $\eta_{A_1} +
\eta_{A_2}$.  Therefore, the change in the process $2Z_{h/2}(t) -
Z_h(t)$ would be $\nu_k$ times
\begin{equation*}
  \eta_{A_1} + \eta_{A_2} + 2\eta_{A_3}.
\end{equation*}
On the other hand, the change in the process generated by the \algNam
Algorithm due to this $k$ is $\nu_k$ times
\begin{equation*}
  \eta_{A_1} + \eta_{A_2} + \eta_{A_3} + \eta_{A_4}.
\end{equation*}
Therefore, and as expected, the means should be the same, but the
variances should not as
\begin{equation*}
  Var ( 2\eta_{A_3}) = 4 Var(\eta_{A_3}) = 2 Var(\eta_{A_3} + \eta_{A_4}).
\end{equation*}
Similarly, in the case in which $\sigma_k^2(X(t))$ decreases as depicted in
\eqref{fig:RE}(b), the process $2Z_{h/2}(t) - Z_h(t)$ would use
$\eta_{A_1} + \eta_{A_2} + \eta_{A_3} - \eta_{A_4}$, whereas the
\algNam Algorithm would use $\eta_{A_1} + \eta_{A_2}$.  Again, the means will be the same, but the variances will not.   In both cases,
the \algNam Algorithm makes better use of the areas to approximate the
quadratic variation of the true process, and thus achieves a higher
order of convergence.

\section{Extension to General Uniformly Elliptic Systems}\label{uniformE}
For a moment let us consider the setting of  general uniformly elliptic SDEs
\begin{equation}
 \begin{aligned}
 dX(t) &= b(X(t))dt + \sum_{k = 1}^M g_k(X(t)) \ dW_k(t),\\
 X(0)&= x \in \R^d
 \end{aligned}
  \label{eq:mainElip}
\end{equation}
where $b$ and $W$ are as before and $g_k\colon\R^d \to \R^d$ is such that
if $G(x)=(g_1(x),\cdots,g_M(x)) (g_1(x),\cdots,g_M(x))^T$ then there
exist positive  $\lambda_-$ and $\lambda_+$ such that
\begin{align*}
\lambda_- |\xi|^2 \leq G(x)\xi \cdot \xi \leq \lambda_+ |\xi|^2  
\end{align*}
for all $x, \xi \in \R^d$. For such a family of uniformly elliptic
matrices a lemma of Motzkin and Wasow \cite{MotzkinWasow53}, whose
precise formulation we take form Kurtz \cite{Kurtz80}, states that if
the entries of $G$ are $C^k$ then there exists an $M$ and $\{
\sigma_k\colon \R^d \to \R_{\ge 0} : k =1,\dots, M\}$, $\{ \nu_k \in
\R^d : k =1,\cdots,M\}$ with $\sigma_k \in C^k$ and strictly positive
so that
\begin{align*}
  G(x) = \sum \sigma_k^2(x) \nu_k\nu_k^T \,.
\end{align*}
Hence \eqref{eq:mainElip} has the same law on path space as
\eqref{eq:main} with these $\sigma_k$ and $\nu_k$. Of course $M$ might
be arbitrarily large (depending on the ratio of $\lambda_+
/\lambda_-$) and hence it is more subtle to compare the total work for
our method with a standard scheme based directly on
\eqref{eq:mainElip}. Furthermore, depending on the dependence on $x$,
it is not transparent how to obtain the vectors $\nu$ and functions
$\sigma$ exactly. Approximations could be obtained using the SVN of
the matrix $G(x)$ for fixed $x$ but we do not explore this further here.
\section{Conclusions and Further Extensions}
\label{sec:conclusion}
We have presented a relatively simple method directly applicable to a
wide class of systems which is weakly second order. We have also shown
how, at least theoretically, it should be applicable to systems which
do not satisfy our structural assumptions but are uniformly
elliptic. We have picked a particularly simple setting to perform our
analysis to make the central points clearer. The assumption that $b$
and $\sigma_k$ are uniformly bounded can be relaxed to a local
Lipschitz condition.  That is to say, if $b$ and $\sigma$ and their
needed derivatives are not bounded uniformly, but rather are bounded
by an appropriate Lyapunov function, then it should be possible to
extend the method directly to the setting of unbounded coefficients
provided the method is stable for the given SDE (see for instance
\cite{MattinglyStuartHigham02}). If the SDE is not globally Lipschitz
then using an implicit drift split-step method as in
\cite{MattinglyStuartHigham02}, an adaptive method as in
\cite{LambdaMattinglyStuart07}, or a truncation method as in
\cite{MilsteinTretyakov05} should extend to our current setting. More
interesting is relaxing the non-degeneracy assumption on the
$\sigma_k$, which was used to minimize the probability of the
diffusion correction being negative. This tact is in some ways
reminiscent of \cite{MilsteinTretyakov05} in that a modification of
the method is made on a small set of paths, though the take here is
quite different. It would be interesting to use the probability that
the correction to the diffusion is negative to adapt the step-size
much in the spirit of \cite{LambdaMattinglyStuart07}. Lastly, there is
some similarity of our method with predictor corrector methods.  In
the deterministic setting, predictor corrector methods not only have a
higher order of accuracy but also have better stability properties.
There have been a number of papers exploring this in the stochastic
context (see \cite{BurrageTian02,NicolaEckhard08,Platen95,Higham00}).
It would be interesting to do the same with the method presented here.

\subsubsection*{Acknowledgments} 
DFA was supported through grant NSF-DMS-0553687 and JCM through grants
NSF-DMS-0449910 and NSF-DMS-06-16710 and a Sloan Foundation
Fellowship.  We would like to thank Andrew Stuart for useful comments
on an early draft and Martin Clark for stimulating questions about
Richardson Extrapolation.  We also thank Thomas Kurtz for pointing out that all uniformly elliptic SDEs can be represented in the form considered in this paper.

\appendix

\section{Proof of Lemma~\ref{lem:keyEstimate}}
\label{app:proof}
The proof of Lemma~\ref{lem:keyEstimate} requires the replacement of
the terms of the form $[\alpha_1 \sigma_k^2(y^*) - \alpha_2
\sigma_k^2(x_0)]^+$ with $[\alpha_1 \sigma_k^2(y^*) - \alpha_2
\sigma_k^2(x_0)]$. The following two lemmas show that this can be done
at the cost of an error whose size is $O(h^3)$. Here $O(h^3)$ has the
same meaning described earlier around \eqref{eq:Onotation}.  We begin
with an abstract technical lemma where $p$ and $q$ satisfy $1/p + 1/q
= 1$.

\begin{proposition}\label{prop:keepPos}
  Let $X$ and $Y$ be a real valued random variables on a probability
  space $(\Omega,\PP)$ with $|XY|_{L^p(\Omega)} < \infty$ for some $p
  \in (1, \infty]$. Then $| \E Y[X]^+ -\E YX| \leq |YX|_{L^p(\Omega)}
  (\PP\{ X < 0\})^{1/q}$. Similarly if $X$,$Y$ and $Z$ are real valued
  random variables with $|ZXY|_{L^p(\Omega)} < \infty$ and $A=\{ X<0\}
  \cup\{Z <0\}$ then $| \E Y[X]^+[Z]^+ -\E YXZ| \leq 2
  |ZYX|_{L^p(\Omega)} (\PP\{A\})^{1/q}$.
\end{proposition}

\begin{proof}
  Let $A=\{ X < 0\}$ and $q=p/(p-1)$.  Then $|\E Y( [X]^+ - X)| \leq
  \E |Y||[X]^+ -X|\one_A \leq |YX|_{L^p(\Omega)} (\PP(A))^{1/q}$,
  showing the first claim. For the second notice that $\E Y[X]^+[Z]^+
  - \E YXZ = (\E Y[X]^+Z - \E YXZ) + (\E Y[X]^+[Z]^+ - \E Y[X]^+Z)$
  and that each of the terms in parentheses can be bounded by the
  first result.
\end{proof}
\begin{corollary}\label{cor:absSwitch} Let $\sigma_k \in \C^2$ with
  $\inf_x \sigma_k(x) >0$ for all $k$ and let $Y$ be a random variable
  with $|Y| \leq C$ a.s. for some $C$. Then for any $p\geq 1 $ there
  exists an $h_0$ so that
  \begin{align*}
    \E Y [\alpha_1 \sigma_k^2(y^*) -& \alpha_2 \sigma_k^2(x_0)]^+ =\E
    Y [\alpha_1 \sigma_k^2(y^*) - \alpha_2 \sigma_k^2(x_0)] + O(h^p)\\
    \E Y [\alpha_1 \sigma_k^2(y^*) -& \alpha_2 \sigma_k^2(x_0)]^+
    [\alpha_1 \sigma_\ell^2(y^*) - \alpha_2 \sigma_\ell^2(x_0)]^+ =\E
    Y [\alpha_1 \sigma_k^2(y^*) - \alpha_2 \sigma_k^2(x_0)] [\alpha_1
    \sigma_\ell^2(y^*) - \alpha_2 \sigma_\ell^2(x_0)] + O(h^p)
  \end{align*}
  for all $h \in (0,h_0]$ and $k,\ell \in \{1,\dots,M\}$, where $y^*$
  is defined via Step 1 of  the \algNam Algorithm.
\end{corollary}
\begin{proof}
  Define the event $\displaystyle A_k = \{\sigma_k(y^*) <
  \frac{\alpha_2}{\alpha_1} \sigma_k(x_0)\}$. In light of
  Proposition~\ref{prop:keepPos}, it is sufficient to show that for
  any $p > 1$ there exists a $C_{p}$ such that $\PP(A_k) \le C_p h^p$.
  Because $\sigma_k$ is Lipschitz there exists a positive $C$ such
  that
  \begin{equation*}
    \sigma_k^2(x_0 + \delta) - \frac{\alpha_2}{\alpha_1}
    \sigma_k^2(x_0) > (1 - \frac{\alpha_2}{\alpha_1}) \sigma_k^2(x_0)
    - C|\delta|, 
  \end{equation*}
  for any $\delta >0$.  In particular, setting $\delta = y^* - x_0 =
  b(x_0)\theta h + \sum_j \sigma_j(x_0)\sqrt{\theta h}\, \nu_j
  \,\eta_{1j}^{(1)}$, and noting that $\alpha_2 < \alpha_1$ and that
  the $\sigma$'s are uniformly bounded from both above and below, the
  result follows from the Gaussian tails of the $\eta$'s.
 \end{proof}

\begin{proof}(of Lemma~\ref{lem:keyEstimate})
  From Taylor's theorem and the definition of the operators involved
  one has
  \begin{align*}
    \E f(\ym) &= f(x_0) + (B_1f)(x_0)\theta h + (B_1^2
    f)(x_0) \frac{\theta^2 h^2}{2} + O(h^3)\\
    &= f(x_0) + (Af)(x_0) \theta h + (B_1^2f)(x_0) \frac{\theta ^2 h^2}{2} +
    O(h^3)\,.
  \end{align*}
  In the last line, we have used the observation that $(B_1f)(x_0) =
  (Af)(x_0)$.  Now we turn to $\E (B_2f)(\ym)$.  We begin by using
  Lemma~\ref{cor:absSwitch} to remove the $[\ccdot]^+$.  Then we use
  the fact that $\alpha_1 - \alpha_2 = 1$ and Taylor's theorem to
  expand various terms to produce the following:
  \begin{align*}
    \E (B_2f)(\ym)&= \E f'(\ym)[\alpha_1 b(\ym) - \alpha_2 b(x_0)] +
    \frac{1}{2} \E\sum_k [\alpha_1 \sigma_k^2(\ym) - \alpha_2
    \sigma_k^2(x_0)]^+
    f''[\nu_k, \nu_k](\ym)\\
    &= \E f'(\ym)[\alpha_1 b(\ym) - \alpha_2 b(x_0)] + \frac{1}{2}
    \E\sum_k [\alpha_1 \sigma_k^2(\ym) - \alpha_2 \sigma_k^2(x_0)]
    f''[\nu_k, \nu_k](\ym)+O(h^2)\\
    &= f'(x_0)[b(x_0)] + \frac{1}{2}  \sum_k \sigma_k(x_0)^2
    f''(x_0)[\nu_k, \nu_k] \\
    & \hspace{.2in} + \E B_1 \Big( f'[\alpha_1 b - \alpha_2 b(x_0)] +
      \frac{1}{2} \sum_k (\alpha_1 \sigma_k^2 - \alpha_2
      \sigma_k^2(x_0)) f''[\nu_k, \nu_k] \Big)
    (x_0) \theta h + O(h^2) \notag \\
    &= (Af)(x_0) + \alpha_1 (B_1(Af))(x_0)\theta h - \alpha_2 (B_1^2
    f)(x_0)\theta h + O(h^2)\notag \\ 
    &= (Af)(x_0)+ \alpha_1 (A^2f)(x_0)\theta h -
    \alpha_2(B_1^2f)(x_0)\theta h + O(h^2).
  \end{align*}
  %
  Similar reasoning produces
  \begin{align*}
    \E (B_2^2 f)(\ym)&= \E \Big( B_2 \big( f'[\alpha_1 b(\ym) -
        \alpha_2 b(x_0)] + \frac{1}{2} \sum_k [\alpha_1
        \sigma_k^2(\ym) - \alpha_2 \sigma_k^2(x_0)]^+
        f''[\nu_k, \nu_k] \big)(\ym)  \Big)\\
    &= f''[b(x_0),b(x_0)](x_0) + \E \sum_k [\alpha_1 \sigma_k^2(\ym) -
    \alpha_2 \sigma_k^2(x_0)]^+ f'''[\nu_k, \nu_k,b(x_0)](x_0) \\
    &+ \frac{1}{4} \E \sum_{k,j} [\alpha_1 \sigma_k^2(\ym) - \alpha_2
    \sigma_k^2(x_0)]^+ [\alpha_1 \sigma_j^2(\ym) -
    \alpha_2 \sigma_j^2(x_0)]^+f''''[\nu_k,\nu_k,\nu_j,\nu_j](x_0) + O(h)\\
    &= f''[b(x_0),b(x_0)](x_0) + \sum_k \sigma_k^2(x_0) f'''[\nu_k,
    \nu_k,b(x_0)](x_0) \\ 
    & \qquad+ \frac{1}{4}
    \sum_{k,j} \sigma_k^2(x_0)  \sigma_j^2(x_0)
    f''''[\nu_k,\nu_k,\nu_j,\nu_j](x_0) + O(h)\\ 
    &= (B_1^2f)(x_0) + O(h)\,.
  \end{align*}
  Combining these estimate and the fact that $2(1-\theta)\theta
  \alpha_2= \theta^2 + (1-\theta)^2$ and $2(1-\theta)\theta
  \alpha_1=1$, produces the quoted result after some algebra.
\end{proof}

\section{Operator Bound for $\mathcal{P}_t\colon \C^k \rightarrow \C^k$}
\label{opBound}

In this section, we show that if $b, \sigma_\ell \in C^k$ then
$\mathcal{P}_t$ is a bounded operator from $\C^m$ to $\C^m$ for $m\in
\{0,\cdots,k\}$. The $k=0$ case follows immediately from $|f(x)| \leq
\|f\|_0$ for all $x \in \R^d$. To address the higher $k$, we introduce
the first $k$ variations of equation \eqref{eq:main}.

For any $\xi \in \R^d$ we denote the first variation of
\eqref{eq:main} in the direction $\xi$ by $J^{(1)}(t,x)[\xi]$ which
solves the linear equation
\begin{align*}
 dJ^{(1)}(t,x)[\xi]&= (\nabla b)(X(t))[J^{(1)}(t,x)[\xi]]\, dt +
\sum_{k=1}^M\nu_k(\nabla \sigma_k)(X(t))[J^{(1)}(t,x)[\xi]]\,
dW_k(t)\,,\\
J^{(1)}(0,x)[\xi]&=\xi \quad\text{and}\quad X(0)=x
\end{align*}
Similarly for $\xi=(\xi_1,\xi_2)
\in\R^2$ the second variation of $X(t)$ (in the directions $\xi$) will
be denoted by $J^{(2)}(t,x)[\xi]$ and defined by 
\begin{align*}
  d J^{(2)}(t,x)[\xi]&= (\nabla b)(X(t))[J^{(2)}(t,x)[\xi]]\, dt +
\sum_{k=1}^M\nu_k(\nabla \sigma_k)(X(t))[J^{(2)}(t,x)[\xi]]\,
dW_k(t) \\ 
& +  (\nabla^2b)(X(t))[J^{(1)}(t,x)[\xi_1],J^{(1)}(t,x)[\xi_2]] +\sum_{k=1}^M(\nabla^2\sigma_k)(X(t))[J^{(1)}(t,x)[\xi_1],J^{(1)}(t,x)[\xi_2]]dW_k(t)\\
J^{(2)}(0,x)[\xi]&=0 \quad\text{and}\quad X(0)=x\,.
\end{align*}
These equations were obtained from successive formal differentiation of
\eqref{eq:main}. By further formal differentiation  we obtain
analogous equations for the $k$-variation $J^{(k)}(t,x)[\xi]$ where
$\xi=(\xi_1,\cdots,\xi_k) \in \R^k$ is the vector of directions.
It is a standard fact that if the coefficients $b, \sigma_j$ are in $\C^k$
then for any $t > 0$
\begin{align*}
\sup_x\E_x \sup \big\{ \sup_{s \in[0,t]} |J^{(n)}(s,x)[\xi_1,\dots,\xi_n]|^p :\xi_i
 \in \R^d \text{ with } |\xi_i|=1\big\} < \infty \,.
\end{align*}
This can be found in Lemma~2 in \cite{Bell87} on p. 196 or in a
slightly different context in Proposition~1.3 in
\cite{Norris86}\footnote{The statement of the Proposition demands
  coefficients in $\C^\infty$. However the bounds on $J^{(k)}$
  only require $\C^k$ coefficients.}. 
With these definitions in hand, we have that for any $f \in C^1$ that
\begin{align*}
\nabla (\mathcal{P}_tf)(x)[\xi]&=\E_x f'(X(t))[J^{(1)}(t,x)[\xi]]\,,\\
\nabla^2 (\mathcal{P}_tf)(x)[\xi]&=\E_x f'(X(t))[J^{(2)}(t,x)[\xi]]+\E_x f^{(2)}(X(t))[J^{(1)}(t,x)[\xi_1],J^{(1)}(t,x)[\xi_2]]\,.
\end{align*}
Using the moment bounds we have that for $q \geq 1$ and an
ever changing constant $C$, 
\begin{align*}
\E \sup_{|\xi|=1}|\nabla (\mathcal{P}_tf)(x)[\xi]|^q \leq&
C\|f\|_{\C^1}^q  \sup_{|\xi|=1}\big|J^{(1)}_t[\xi]\big|^q \leq C\|f\|_{\C^1}^q < \infty \\
\E \sup_{|\xi_i|=1}|\nabla^2 (\mathcal{P}_tf)(x)[\xi_1,\xi_2]|^q \leq&
C\|f\|_{\C^2}^q\Big(\big(\E \sup_{|\xi_1|=1}
|J^{(1)}(t,x)[\xi_1]|^{2q}\big)^{\frac12}+ \E \sup_{|\xi_i|=1}
|J^{(2)}(t,x)[\xi_1,\xi_2]]|^q\Big)\\
\leq& C\|f\|_{\C^2}^q < \infty
\end{align*}
Continuing in this manner we see that for any positive integer $m$ if
$f,b, \sigma_\ell \in \C^m$ then for any $q \geq
1$  one has 
\begin{align*}
\E \sup_{|\xi_i|=1}|\nabla^m (\mathcal{P}_tf)(x)[\xi_1,\cdots,\xi_m]|^q
\leq C\|f\|_{\C^m}^q < \infty\,
\end{align*}
for some $C$.  
Now observe that taking $q=1$ proves the desired claim on the operator
norm of $\mathcal{P}_t$ from $\C^k$ to $\C^k$ since
\begin{align*}
\|\mathcal{P}_t f\|_{k}& \leq C \sum_{j=0}^k  \E
\sup_{|\xi_i|=1}\big|(\nabla^j \mathcal{P}_tf)(x)[\xi_1,\cdots,\xi_j]\big|
\leq C \sum_{j=0}^k \|f\|_{\C^j} \leq C\|f\|_{\C^k}
\end{align*}\,.

\bibliographystyle{siam} 

\def\cprime{$'$} \def\cprime{$'$}

  \end{document}